\documentclass[a4paper]{article}
\usepackage[latin1]{inputenc}
\usepackage[english]{babel}
\usepackage{amsmath}
\usepackage{amssymb,amsthm,enumerate}

\usepackage[bookmarks,bookmarksopen,colorlinks]{hyperref}
\newcommand{\R}{\mathbb{R}}
\newcommand{\C}{\mathbb{C}}
\newcommand{\N}{\mathbb{N}}
\newcommand{\Z}{\mathbb{Z}}

\newcommand{\im}{\operatorname{Im}}
\newcommand{\indic}{\operatorname{l\negthinspace l}}

\newcommand{\ja}[2]{j_{#1}\left(\lambda  {#2}\right)}
\newcommand{\na}[2]{\eta_{#1}\left(\lambda  {#2}\right)}
\newcommand{\n}[2]{{\left\|{#1}\right\|}_{#2}}
\newcommand{\ing}{[\negthinspace [}
\newcommand{\ind}{]\negthinspace ]}
\newcommand{\vect}[1]{\operatorname{Span}\left\{#1\right\}}

\newcommand{\G}{\mathcal{G}}
\newcommand{\classe}[1]{\mathcal{C}^{#1}}
\newcommand{\lc}{L_{\C}^2(0,1)}
\newcommand{\lr}{L_{\R}^2(0,1)}

    \newtheorem{thm}{Theorem}[section]
    \newtheorem{coro}{Corollary}[section]
    \newtheorem{prop}{Proposition}[section]
    \newtheorem{lem}{Lemma}[section]

\theoremstyle{definition}
    \newtheorem{defn}{Definition}[section]
    \newtheorem*{nota}{Notations}
\theoremstyle{remark}
    \newtheorem{rem}{Remark}

    \numberwithin{equation}{section}
    \newcommand{\be}{\begin{equation*}}
    \newcommand{\ee}{\end{equation*}}
\renewcommand{\wr}{\mathcal{W}}
\renewcommand{\O}[1]{\mathcal{O}\negthinspace\left(#1\right)}

\author{Frédéric SERIER\\
\small\emph{Laboratoire de Math\'{e}matiques Jean Leray}\\
\small\emph{ UMR CNRS-Universit\'{e} de Nantes,%}\\
Facult\'{e} des sciences et techniques,}\\
\small\emph{2 rue de la Houssini\`{e}re, BP 92208, 44322 Nantes
cedex 03, France}}
\title{Inverse spectral problem for radial Schrödinger operator on $[0,1]$.\\}

\begin{document}
    \maketitle
        \begin{abstract}
For a class of singular Sturm-Liouville equations on the unit
interval with explicit singularity $a(a + 1)/x^2, a \in \N$, we
consider an inverse spectral problem. Our goal is the global
parametrization of potentials by spectral data noted by
$\lambda^a$, and some norming constants noted by $\kappa^a$. For
$a = 0$ and $1,$ $\lambda^a\times \kappa^a$ was already known to
be a global coordinate system on $\lr$. With the help of
transformation operators, we extend this result to any
non-negative integer $a$ and give a description of isospectral
sets.
\end{abstract}

\section{Introduction}

Inverse spectral problem for Schrödinger operator with radial
potential $\mathrm{H}:=-\Delta+q(\n{X}{})$ acting on the unit ball
of $\R^3$ leads by separation of variables (see \cite{rs}, p.
$160-161$) to consider a collection of singular differential
operators $\mathrm{H}_a(q)$, $a\in \N$ acting on $\lr$, defined by

\begin{equation}\label{SL+V}
    \mathrm{H}_a(q) y(x):=\left(-\frac{d^2}{d x^2}+\frac{a(a+1)}{x^2}+q(x)
    \right)y(x)=\lambda y(x),\quad x\in[0,1], \lambda\in\C,
\end{equation}
with Dirichlet boundary conditions
\begin{equation}\label{SL-bord}
    y(0)=y(1)=0.
\end{equation}

Our purpose is the stability of the inverse spectral problem for
$\mathrm{H}_a$. For this, we construct for each $a\in\N$, a
standard map $\lambda^a\times\kappa^a$ for potentials $q\in\lr$
with spectral data $\lambda^a$ and some terminal velocities
$\kappa^a$.

Borg \cite{bo} and Levinson \cite{lev} proved that
$\lambda^0\times\kappa^0$ is one-to-one on $\lr$. Pöschel and
Trubowitz \cite{ist} improved this result obtaining
$\lambda^0\times\kappa^0$ as a global real-analytic coordinate
system on $\lr$. Guillot and Ralston \cite{gr} extended their
results to $\lambda^1\times\kappa^1$. Next Carlson \cite{carlson}
(see also Zhornitskaya and Serov \cite{zs}) proved that for all
real $a\geq -1/2$, $\lambda^a\times\kappa^a$ is one-to-one on
$\lr$.

We complete these works, to any integer $a$, with the realization
of $\lambda^a\times\kappa^a$ as a local real-analytic coordinate
system on $\lr$. Moreover, we obtain a description for the set of
potentials with same spectral data : these so-called isospectral
sets are proved to be regular manifolds, expression for their
tangent and normal spaces are given. The key point of the proof
lies in the use of transformation operators: they help us to
handle Bessel functions naturally underlying in this problem.

\section{The direct spectral problem.}

Part of the following properties are deduced or came from
\cite{gr}, \cite{carlson}-\cite{carlson2} and \cite{zs};
nevertheless main structure is given by \cite{ist}.

Define $\omega:=\lambda^{1/2}$ (the square root determination is
pointless because of the parity of functions using it) and for all
$n\in\N$
$$\displaystyle(2n+1)!!:=1\cdot 3\cdots
(2n+1).$$

A fundamental system of solutions for \eqref{SL+V} when $q=0$ is
given by
$$u(x,\lambda)=\frac{(2a+1)!!}{\omega^{a+1}}\,j_a(\omega x),\quad\displaystyle
v(x,\lambda)=-\frac{\omega^{a}}{(2a+1)!!}\, \eta_a(\omega x), $$
where $j_a$ and $\eta_a$ are spherical Bessel functions (see
section \ref{annexe-bessel}). This family is called fundamental
since its wronskian $\wr(u,v)$ is equal to $1$. From their
behavior near $x=0$, $u(x,\lambda)$ is called the regular
solution, it is analytic on $[0,1]\times\C$; $v(x,\lambda)$ is
called the singular solution, it is analytic on
$(0,1]\times\C$. %(cf. section \ref{annexe-bessel-series})

Following Guillot and Ralston \cite{gr} but also Zhornitskaya and
Serov \cite{zs}, we construct solutions for \eqref{SL+V} by a
Picard's iteration method from $u$ and $v$.

Let $\varphi$ and $\widetilde{\psi}$ be defined by
$$\varphi(x,\lambda,q)=\sum_{k\geq 0}\varphi_k(x,\lambda,q),\quad\widetilde{\psi}(x,\lambda,q)=\sum_{k\geq 0}\widetilde{\psi}_k(x,\lambda,q)$$
with
\begin{gather}\label{sl-relation_recurrence_reg}\left\{\begin{array}{c}
   \varphi_0(x,\lambda,q)=u(x,\lambda), \\
      \varphi_{k+1}(x,\lambda,q)=\displaystyle\int_0^x
    \G(x,t,\lambda)q(t)\varphi_k(t,\lambda,q)dt,\quad k\in\N;
\end{array}\right.\\\label{sl-relation_recurrence_sing}
\left\{\begin{array}{c}
    \widetilde{\psi}_0(x,\lambda,V)=v(x,\lambda), \\
    \widetilde{\psi}_{k+1}(x,\lambda,q)=\displaystyle-\int_x^1
    \G(x,t,\lambda)q(t)\widetilde{\psi}_k(t,\lambda,q)dt,\quad k\in\N.
\end{array}\right.
\end{gather}
$\G$ is called Green function and is given by
\begin{equation*}%\label{sl-resolvante}
    \G(x,t,\lambda)=v(x,\lambda)u(t,\lambda)-u(x,\lambda)v(t,\lambda), \quad (x,t)\in(0,1)\times(0,1).
\end{equation*}

Now follows the expected result:
%%%%%%%%%%%%%%%%%%%%%%%%%%%%%%%%%%%%%%%%%%%%%%%%%%%%%%%%%%%%%%%%%%%%%%%%%%%
\begin{lem}%\label{sl-lemme_base}
    Series defined in \eqref{sl-relation_recurrence_reg}, respectively in \eqref{sl-relation_recurrence_sing} uniformly converge
    on bounded sets of $[0,1]\times\C\times\lc$, respectively of
    $(0,1]\times\C\times\lc$ towards solutions of \eqref{SL+V}.
    Moreover, they satisfy the integral equations
    \begin{gather}
        \label{sl-eq-int-phi}\varphi(x,\lambda,q)=u(x,\lambda)+\int_0^x
            \G(x,t,\lambda)q(t)\varphi(t,\lambda,q)dt,\\
        \label{sl-eq-int-psitilde}\widetilde{\psi}(x,\lambda,q)=v(x,\lambda)-\int_x^1
            \G(x,t,\lambda)q(t)\widetilde{\psi}(t,\lambda,q)dt,
    \end{gather}
    and the estimates
    \begin{eqnarray*}
      \left|\varphi(x,\lambda,q)\right| &\leq& C e^{|\im
            \omega|x}\left(\frac{x}{1+|\omega|x}\right)^{a+1} , \\
      \left|\widetilde{\psi}(x,\lambda,q)\right|  &\leq& C e^{|\im
            \omega|(1-x)}\left(\frac{1+|\omega|x}{x}\right)^a,
    \end{eqnarray*}
with $C$ uniform on bounded set of $\lc$.
\end{lem}
%%%%%%%%%%%%%%%%%%%%%%%%%%%%%%%%%%%%%%%%%%%%%%%%%%%%%%%%%%%%%%%%%%%%%%%%%%%
\begin{proof}
The proof is given for $\varphi$, it is the same for
$\widetilde{\psi}$. Estimates \eqref{annexe-bessel_est_ja} for
Bessel functions gives
\begin{equation}\label{sl-estimation-R}
    |u(x,\lambda)|\leq C e^{|\im\omega|x}
    \left(\frac{x}{1+|\omega|x}\right)^{a+1}.
\end{equation}
Iterative relation \eqref{sl-relation_recurrence_reg} leads to
\begin{equation*}%\label{sl-exp_R1}
    \varphi_1(x,\lambda,q)=\int_0^x \G(x,t,\lambda)q(t)u(t,\lambda)dt,
\end{equation*}
which, combining \eqref{sl-estimation-R} and the Green function
estimate \eqref{annexe-resolvante-sl_est_tx}, is bounded by
\begin{equation*}
    \left|\varphi_1(x,\lambda,q)\right|\leq C^2 e^{|\im\omega|x}\left(\frac{x}{1+|\omega|x}\right)^{a+1} \int_0^x \frac{t|q(t)|}{1+|\omega|t}dt,
\end{equation*}
By successive iterations and recurrence, for all positive integer
$n$, we get
\begin{equation*}
    \left|\varphi_n(x,\lambda,q)\right|\leq \frac{C^{n+1}}{n!}
    e^{|\im\omega|x}\left(\frac{x}{1+|\omega|x}\right)^{a+1}
     \left(\int_0^{x} \frac{t|q(t)|}{1+|\omega|t}dt\right)^n.
\end{equation*}
This proves uniform convergence on bounded sets of
$[0,1]\times\C\times\lc$ for $\varphi$ and the estimate. Integral
relation follows from \eqref{sl-relation_recurrence_reg}.
\end{proof}

\begin{coro}
$\varphi$ and $\widetilde{\psi}$, solutions for \eqref{SL+V}
follow the estimates
    \begin{gather}
        \label{sl-estimation-R-L2}\left|\varphi(x,\lambda,q)-u(x,\lambda)\right| \leq C e^{|\im
            \omega|x}\left(\frac{x}{1+|\omega|x}\right)^{a+1} \int_0^x \frac{t|q(t)|}{1+|\omega|t}d t,\\
        \label{sl-estimation-S-L2}\left|\widetilde{\psi}(x,\lambda,q)-v(x,\lambda)\right| \leq C e^{|\im
            \omega|(1-x)}\left(\frac{1+|\omega|x}{x}\right)^a \int_x^1 \frac{t|q(t)|}{1+|\omega|t}d
            t.
    \end{gather}
    Let $'$ denote the derivative with respect to $x$, we also obtain
    \begin{gather*}
        \varphi'(x,\lambda,q)=u'(x,\lambda)+\int_0^x
            \frac{\partial\G}{\partial x}(x,t,\lambda)q(t)\varphi(t,\lambda,q)dt,\\
        \widetilde{\psi}'(x,\lambda,q)=v'(x,\lambda)-\int_x^1
            \frac{\partial\G}{\partial x}(x,t,\lambda)q(t)\widetilde{\psi}(t,\lambda,q)dt,
    \end{gather*}
    and
    \begin{gather}
        \label{sl-estimation-Rprime}\left|\varphi'(x,\lambda,q)-u'(x,\lambda)\right| \leq C e^{|\im
            \omega|x}\left(\frac{x}{1+|\omega|x}\right)^{a} \int_0^x \frac{t|q(t)|}{1+|\omega|t}d t,\\
        %\label{sl-estimation-Sprime}
        \left|{\widetilde{\psi}}'(x,\lambda,q)-v'(x,\lambda)\right|\leq C e^{|\im
            \omega|(1-x)}\negthinspace\left(\frac{1+|\omega|x}{x}\right)^{a+1}\negthinspace\negthinspace \int_x^1 \frac{t|q(t)|}{1+|\omega|t}d t,
    \end{gather}
    with $C$ uniform on bounded sets of $\lc$.
\end{coro}
\begin{proof}
    For the first estimates, we just have to add those from the proof of the preceding lemma beginning at $k=1$. Integral equations
    follow by derivation and the last estimates follow using estimates \eqref{annexe-diffres-sl_est_tx} and
    \eqref{annexe-diffres-sl_est_xt} for $\G$ given in annexe \ref{annexe-bessel}.
\end{proof}

\begin{rem} Notice that for $\omega\neq 0$,
$$0\leq\frac{t}{1+|\omega|t}=\frac{1}{|\omega|}\frac{|\omega|t}{1+|\omega|t}\leq \frac{1}{|\omega|}.$$
Thus, bound from the corollary are asymptotic estimates when
$\omega\rightarrow\infty$.
\end{rem}
\begin{rem} Theses estimates lead to
\begin{equation}\label{sl-def_sol_sing}
        \wr(\lambda,q):=\wr(\varphi(\cdot,\lambda,q),\widetilde{\psi}(\cdot,\lambda,q))=1+\O{\frac{1}{\omega}},\quad \omega\rightarrow
        \infty.
\end{equation}
\end{rem}

According to \cite{new} and \cite{carlson2}, following \cite{gr},
there exists a function $\psi(x,\lambda,q)$ such that
$x^a\psi(x,\lambda,q)$ is analytic on $\C\times\lc$ for all
$x\in[0,1]$ and such that $\{\varphi,\psi\}$ is a basis for the
solutions of \eqref{SL+V}, with the normalization
$$\wr(\varphi,\psi)=1.$$

Behavior near $x=0$ for $\varphi$ is thus inherited from $u$ and
with the wronskian relation this leads to the following boundary
conditions at $x=0$
  \begin{gather*}
    %\label{sl-Condition-en-0}
    \lim_{x\rightarrow 0^+}
    \frac{\varphi(x,\lambda,q)}{x^{a+1}}=1,\quad
    \lim_{x\rightarrow 0^+}
    x^a \psi(x,\lambda,q)=1.
\end{gather*}
Thus, $\psi$ is called as the singular solution of
\eqref{operateur_reduction}.

With the estimate \eqref{sl-def_sol_sing}, $\wr(\lambda,q)$ has no
zeros near the real axis and we can take for this singular
solution the one defined by
\begin{equation}\label{def-psi-bis}
    \psi(x,\lambda,q)=\frac{\widetilde{\psi}(x,\lambda,q)}{\wr(\lambda,q)}.
\end{equation}

For real values of $a$, it is interesting to read \cite{carlson}
and \cite{zs} (precisely when $a\geq\frac{-1}{2}$), particularly
non integer values allow to study the spectral problem for the
radial Schrödinger operator acting on the unit ball in a even
dimension space (in our case, the operator acts on $\R^3$).

Uniform convergence for series in the previous lemma give
regularity for $\varphi$ and $\tilde{\psi}$, expressed as bellow

\begin{prop}[Analyticity of solutions]$\quad$
\begin{enumerate}[(a)]
    \item  For all $x\in[0,1]$, $\varphi(x,\lambda,q)$
        is analytic on $\C\times\lc$. Moreover, it is real valued on $\R\times\lr$.
    \item The map $$\varphi: (\lambda,q)\mapsto \varphi(\cdot,\lambda,q)$$
    is analytic from $\C\times\lc$ to
    $H^2([0,1],\C)$.
    \item For all $x\in(0,1]$, $\tilde{\psi}(x,\lambda,q)$ is analytic on
    $\C\times\lc$ and real valued on $\R\times\lr$.
\end{enumerate}
\end{prop}

Consequently, $\wr(\lambda,q)$ is an analytic function on
$\C\times\lc$ and thus $\psi$ defined by \eqref{def-psi-bis} has
the needed regularity.

Regularity leads to existence of derivatives, precisely we obtain
expression for $\lc$-gradients with respect to $\lambda$ and $q$.
For this, following \cite{ist}, we calculate differential for
$\varphi$ with respect to $q$:
%   Pour cela, rappelons brièvement la notion de gradient dans un
%   espace de Hilbert:
%   \begin{rem}[Rappel] Si $E$ est un espace de Hilbert et $F$ est la droite réelle ou complexe, et si $f:E\rightarrow
%   F$ est différentiable, le gradient de $f$ relativement à $x\in E$
%   est le représentant de Riesz de la différentielle de $f$ au point
%   $x$, c'est à dire $$d_x f(v)=\left\langle v,\overline{\nabla_x
%   f}\right\rangle_E,\quad\forall v\in E.$$
%   \end{rem}
$$\left[d_{q}
\varphi(x,\lambda,q)\right](v)=\int_0^x
    \tilde{\G}(x,t,\lambda,q) v(t)\varphi(t,\lambda,q) dt,$$
where
$$\tilde{\G}(x,t,\lambda,q)=\psi(x,\lambda,q)\varphi(t,\lambda,q)-\varphi(x,\lambda,q)\psi(t,\lambda,q).$$
Rewriting this relation as a scalar product, we obtain
\begin{prop}\label{sl-gradient-sol}[Gradients for the regular solution]$\,$\\
For all $v\in\lc$, we have
\begin{gather*}
    %\label{sl-differentielle_R_V}
    \nabla_{q}\,
    \varphi(x,\lambda,q)(t)=\varphi(t,\lambda,q)
    \left[\psi(x,\lambda,q)\varphi(t,\lambda,q)-\varphi(x,\lambda,q)\psi(t,\lambda,q)\right]
    \indic_{[0,x]}(t),\\
    %\label{sl-gradient_R_lambda}
    \frac{\partial \varphi}{\partial
    \lambda}(x,\lambda,q)=-\left[d_{q}
    \varphi(x,\lambda,q)\right](1),\\
    %\label{sl-differentielle_R_V_prime}
    \nabla_{q}\,
    \varphi'(x,\lambda,q)(t)=\varphi(t,\lambda,q)
    \left[\psi'(x,\lambda,q)\varphi(t,\lambda,q)-\varphi'(x,\lambda,q)\psi(t,\lambda,q)\right]
    \indic_{[0,x]}(t).
\end{gather*}
\end{prop}
%Rappelons que l'expression de la dérivée de $\varphi$ par rapport
%à $\lambda$ se déduit de celle du gradient par rapport à $q$ en
%écrivant l'identité
%$$\varphi(x,\lambda+\varepsilon,q)=\varphi(x,\lambda,q-\varepsilon).$$

\subsection{Spectral data}

From this point, potential $q$ has real values, that is
$$q\in\lr.$$

\subsubsection{Eigenvalues}

Spectrum localization and normalization for eigenfunctions is made
following \cite{ist} and \cite{gr}. Boundary conditions given by
\eqref{SL-bord} define spectra of \eqref{SL+V}-\eqref{SL-bord} as
the zeros set of the entire function $\lambda\mapsto
\varphi(1,\lambda,q)$. Results from Zhornitskaya and Serov
\cite{zs} give simplicity and localization similarly to \cite{ist}
and \cite{gr}. Regularity and expression for the gradient for each
eigenvalue follows from \cite{carlson2}.

\begin{nota} Let $\left(\lambda_{a,n}(q)\right)_{n\geq 1}$ be the set of eigenvalues representing the spectrum $\mathrm{H}_a(q)$.
For every integer $n\geq 1$, let $g_n(t,q)$ be the eigenvector
with respect to $\lambda_{a,n}(q)$ defined by the normalization
$$g_n(t,q)=\frac{\varphi(t,\lambda_{a,n}(q),q)}{\n{\varphi(\cdot,\lambda_{a,n}(q),q)}{\lr}}.$$
\end{nota}

Now recall previous results with the following :
\begin{thm} Let $q\in\lr$, the spectrum for the problem \eqref{SL+V} with
    Dirichlet boundary conditions \eqref{SL-bord} is a strictly increasing sequence of eigenvalues
    $\lambda^a(q)=(\lambda_{a,n}(q))_{n\geq 1}$ which are all real-analytic on $\lr$ and verify
    \begin{equation}\label{asymptotique_spectre-bidon}
        \lambda_{a,n}(q)=\left(n+\frac{a}{2}\right)^2\pi^2\left(1+\O{\frac{1}{n}}\right),
    \end{equation}
    \begin{equation*}%\label{gradient_spectre}
        \nabla_q \lambda_{a,n}(t)=g_n(t,q)^2.
    \end{equation*}
\end{thm}

\begin{coro}$\,$\\
Let $q\in\lr$, the following asymptotics
    \begin{equation}\label{sl-est-vect-propre}
        g_n(t,q)=\sqrt{2} j_{a}{(\omega_{a,n}(q)\,
        t)}+\O{\frac{1}{n}},
    \end{equation}
    \begin{equation}\label{gradient_spectre_estimation}
        \nabla_q \lambda_{a,n}(t)=2 j_{a}{(\omega_{a,n}(q)\,
        t)}^2+\O{\frac{1}{n}},
    \end{equation}
    \begin{equation}\label{sl-asymptotique_spectre}
        \lambda_{a,n}(q)=\left(n+\frac{a}{2}\right)^2\pi^2+\O{1},
    \end{equation}
are uniform on bounded sets in $\lr$.
\end{coro}

\begin{proof}
Rewriting relation \eqref{sl-estimation-R-L2} gives for
$\omega\in\R$,
\begin{equation}\label{sl-preuve-est-L2-1}
    \varphi(x,\lambda,q)=\frac{(2a+1)!!}{\omega^{a+1}}j_a(\omega x)
    +\O{\frac{1}{\omega}\left(\frac{x}{1+|\omega|x}\right)^{a+1}}.
\end{equation}
Thus,
$$\int_0^1 \varphi(t,\lambda,q)^2 d t=\frac{[(2a+1)!!]^2}{\omega^{2a+2}}\int_0^1 \left[j_a(\omega
    t) +\O{\frac{1}{\omega}\left(\frac{|\omega|x}{1+|\omega|x}\right)^{a+1}}\right]^2 dt.$$
Estimate \eqref{annexe-bessel_est_ja} for Bessel function $j_a$
leads to
$$\left[j_a(\omega
    t)+\O{\frac{1}{\omega}\left(\frac{|\omega|x}{1+|\omega|x}\right)^{a+1}}\right]^2=j_a(\omega
    t)^2 +\O{\frac{1}{\omega}\left(\frac{|\omega|x}{1+|\omega|x}\right)^{2a+2}},$$
then
 $$\int_0^1 \varphi(t,\lambda,q)^2 d t=\frac{[(2a+1)!!]^2}{\omega^{2a+2}}\left[\int_0^1 j_a(\omega
 t)^2 dt +\O{\frac{1}{\omega}}\right].$$
Relation \eqref{annexe-bessel-int-1} gives
 \begin{equation}\label{sl-preuve-est-L2-2}
    \int_0^1 \varphi(t,\lambda,q)^2 d t=\frac{1}{2}\frac{[(2a+1)!!]^2}{\omega^{2a+2}}\left[1+\O{\frac{1}{\omega}}\right].
\end{equation}
Using both \eqref{sl-preuve-est-L2-1} and
\eqref{sl-preuve-est-L2-2}, knowing that
$\displaystyle\omega_{a,n}(q)=\left(n+\frac{a}{2}\right)\pi+\O{1}$,
relations \eqref{sl-est-vect-propre} and
\eqref{gradient_spectre_estimation} follow.\\
For \eqref{sl-asymptotique_spectre}, we write
\begin{eqnarray*}
  \lambda_{a,n}(q)-\lambda_{a,n}(0) &=& \int_0^1 \frac{d}{dt}\left(\lambda_{a,n}(t q)\right)dt=
        \int_0^1 \left\langle \nabla_{t q}\lambda_{a,n},q\right\rangle dt, \\
   &=&  \iint_{[0,1]^2} 2 j_a{\left(\omega_{a,n}(t q)x\right)}^2 q(x)d x\,
   dt+\O{\frac{1}{n}}.
\end{eqnarray*}
Thus,
\begin{multline}\label{sl-integrale-In}
\lambda_{a,n}(q)=\lambda_{a,n}(0)+\int_0^1 q(t)dt\\
+\iint_{[0,1]^2} \Big[2 j_a{\left(\omega_{a,n}(t
q)x\right)^2}-1\Big]q(x)d x\,   dt+\O{\frac{1}{n}}.
\end{multline}
Eq. \eqref{annexe-bessel_est_ja} shows the boundedness for the
above second integral. According to \cite{rusa} Eq. (2.10), we
have
$$\lambda_{a,n}(0)=\left(n+\frac{a}{2}\right)^2\pi^2-a(a+1)+\ell^2(n),$$
which proves \eqref{sl-asymptotique_spectre}.
\end{proof}

The asymptotic for eigenvalues is far from optimal, but at this
point we cannot do better unless more work to deal with the second
integral in \eqref{sl-integrale-In}.

\subsubsection{Terminal velocities}

Solving inverse spectral problem implies the knowledge of other
spectral data. Indeed even for the regular case, the eigenvalues
don't give the associated potential $q$. Nevertheless, partial
results (as Borg theorem) exist. Thus, let's define the needed
additional data as in \cite{ist} and \cite{gr}.

\begin{defn} Let $\kappa_{a,n}$, called terminal velocity, be defined for all integer $n\geq 1$
by
\begin{equation*}%\label{sl-def-kappa}
    \kappa_{a,n}(q)=\ln\left|\frac{\varphi'(1,\lambda_{a,n}(q),q)}{u'(1,\lambda_{a,n}(0))}\right|.
\end{equation*}
Moreover, let $a_n(t,q)$ be the function given by
    \begin{equation*}%\label{sl-def-an}
    a_n(t,q)=\varphi(t,\lambda_{a,n}(q),q)\psi(t,\lambda_{a,n}(q),q).
    \end{equation*}
\end{defn}

Give some properties for these items:

\begin{thm}
For all $n\in\N,$ $q\mapsto\kappa_{a,n}(q)$ is real-analytic on
$\lr$. Its gradient is expressed by
    \begin{equation}\label{gradient_kappa}
        \nabla_q \kappa_{a,n}(t)=-a_n(t,q)+\nabla_q \lambda_{a,n}(t)\int_0^1 a_n(s,q)d s,
    \end{equation}
and follows the estimate
    \begin{equation}\label{gradient_kappa_estimation}
        \nabla_q \kappa_{a,n}(t)=\frac{1}{\omega_{a,n}} j_{a}{(\omega_{a,n} t)}\eta_a{(\omega_{a,n}
        t)}+\O{\frac{1}{n^2}},
    \end{equation}
uniformly on bounded sets in $[0,1]\times\lr$.
\end{thm}
\begin{proof}[Preuve]Regularity of $\kappa_{a,n}$ comes from regularity of $\lambda_{a,n}$ and $\varphi$.
Expression for its gradient is a straightforward calculation using
Proposition \ref{sl-gradient-sol}.

Determine the asymptotic value of the gradient. Estimations
\eqref{sl-estimation-R-L2}-\eqref{sl-estimation-S-L2} and
\eqref{sl-def_sol_sing} lead to
\begin{multline*}
    a_n(x,q)=\Bigg[u(x,\lambda_{a,n})+\O{\frac{1}{\omega_{a,n}}\left(\frac{x}{1+|\omega_{a,n}x|}\right)^{a+1}}\Bigg]\\
    \times\Bigg[v(x,\lambda_{a,n})+\O{\frac{1}{\omega_{a,n}}\left(\frac{1+|\omega_{a,n}x|}{x}\right)^{a}}\Bigg],
\end{multline*}
then, definitions of $u$ and $v$ give
\begin{multline*}
    a_n(x,q)=\frac{-1}{\omega_{a,n}}\Bigg[j_a(\omega_{a,n}
x)+\O{\frac{1}{\omega_{a,n}}\left[\frac{|\omega_{a,n}|x}{1+|\omega_{a,n}x|}\right]^{a+1}}\Bigg]\\
    \times\Bigg[\eta_a(\omega_{a,n}
x)+\O{\frac{1}{\omega_{a,n}}\left[\frac{1+|\omega_{a,n}x|}{|\omega_{a,n}|x}\right]^{a}}\Bigg].
\end{multline*}
From \eqref{annexe-bessel_est_ja} and
\eqref{annexe-bessel_est_na}, we deduce
$$a_n(t,q)=\frac{-1}{\omega_{a,n}}j_a(\omega_{a,n}x)\eta_a(\omega_{a,n}x)+
\O{\frac{1}{{\omega_{a,n}}^2}}.$$ Now, we have
$$\int_0^1 a_n(t,q)dt=\frac{-1}{\omega_{a,n}}\int_0^1 j_a(\omega_{a,n}t)\eta_a(\omega_{a,n}t)dt+\O{\frac{1}{{\omega_{a,n}}^2}}.$$
Relation \eqref{annexe-bessel-int-2} gives the result.
\end{proof}

As noticed previously, on the contrary to the regular case, good
estimations for terminal velocities are harder to obtain. Indeed,
we have the relation
$$\kappa_{a,n}(q)=\int_0^1 \frac{d}{dt}\left(\kappa_{a,n}(t q)\right) dt
=\int_0^1 \left\langle \nabla_{t q}\kappa_{a,n},
q\right\rangle_{\lr} dt.$$ Using Eq.
\eqref{gradient_kappa_estimation}, we get
\begin{equation}\label{sl-integrale-Jn}
    \kappa_{a,n}(q)=\iint_{[0,1]^2} \frac{j_{a}{\big(\omega_{a,n}(t q)
x\big)}\eta_a\big(\omega_{a,n}(t q)
        x\big)}{\omega_{a,n}(t q)}\,  q(x)d x\, d t+\O{\frac{1}{n^2}}.
\end{equation}
Thus, follows
    \begin{equation*}%\label{sl-est-kappa-bidon}
        \kappa_{a,n}(q)=\O{\frac{1}{n}},
    \end{equation*}
losing the $\O{\frac{1}{n^2}}$ accuracy. As Guillot and Ralston
\cite{gr} for $a=1$, we'll use transformation operator to avoid
this problem.

We finish our study of the direct problem some properties of the
gradients.

\subsection{Orthogonality}

\begin{prop} For all $(n,m)\in \N^2$, $n,m\geq 1$, we have%
\begin{enumerate}
    \item $\displaystyle \left\langle {g_n}^2, \frac{d}{dx}{g_m}^2   \right\rangle=0$,
    \item $\displaystyle \left\langle a_n,\frac{d}{dx}{g_m}^2  \right\rangle=\frac{1}{2}\delta_{n,m}$,
    \item $\displaystyle \left\langle a_n, \frac{d}{dx}{a_m} \right\rangle=0$.
\end{enumerate}
\end{prop}
\begin{proof} Formally, the proof follows as in \cite{ist}, thus we just prove the first one.
Let $(n,m)\in\N^2$, with $n,m\geq 1$, integration by parts gives
$$\left\langle {g_n}^2,
\frac{d}{dx}{g_m}^2\right\rangle=\frac{1}{2}\int_0^1
\left[{g_n}^2\left({g_m}^2\right)'-{g_m}^2\left({g_n}^2\right)'\right]dt=\int_0^1
g_n g_m \wr{(g_n,g_m)}dt. $$%
When $n=m$, the first relation is true. Now suppose $n\neq m$, we
use the identity
$$\frac{d}{d x}\wr(g_n,g_m)=(\lambda_{a,n}-\lambda_{a,m})g_n g_m$$
to obtain
$$\left\langle {g_n}^2,
\frac{d}{dx}{g_m}^2\right\rangle=\frac{1}{2(\lambda_{a,n}-\lambda_{a,m})}\Big[\wr{(g_n,g_m)}^2\Big]_0^1=0.
$$
\end{proof}

Now, we deduce algebraic properties for gradients. For this,
recall the definition of the linear independence in an Hilbert
space.

\begin{defn} An infinite vector family $(u_k)_{k\geq 1}$ in an Hilbert
space is called free or its vectors are linearly independent if
each vector of this family doesn't belong to the closed span of
the others. Precisely
$$\forall k\geq 1\,,\quad u_k\notin \overline{\operatorname{Vect}\left(u_j|j\geq 1, j\neq k\right)}.$$
\end{defn}

%\begin{rem} Cette notion de liberté de vecteurs n'est pas algébrique. Elle est à mettre en parallèle avec
%la notion de base hilbertienne. Rappelons qu'un espace vectoriel
%normé admettant une base dénombrable de vecteurs n'est jamais
%complet.
%\end{rem}

%Nous obtenons suivant \cite{ist} les propriétés suivantes:

\begin{coro}\label{sl-gradients_libres}[Gradients and orthogonality]$\,$
\begin{enumerate}[(a)]
    \item $1$, $
        \left\{{g_n}^2-1\right\}_{n\geq 1}$ are linearly independent, and so are
        $\displaystyle\left\{\frac{d}{d x}{g_n}^2\right\}_{n\geq 1}$. Moreover, these two families are
        orthogonal.
    \item  For all $(n,m)\in\N^2$, $n,m\geq 1$ we have
        \begin{enumerate}[(i)]
        \item $\displaystyle\left\langle\nabla_{q}\kappa_{a,n},\frac{d}{d x}\left(\nabla_{q}\kappa_{a,m}\right)\right\rangle=0,$
        \item $\displaystyle\left\langle\nabla_{q}\kappa_{a,n},\frac{d}{d x}\left(\nabla_{q}\lambda_{a,m}\right)\right\rangle=\frac{1}{2}\delta_{n,m},$
        \item $\displaystyle\left\langle\nabla_{q}\lambda_{a,n},\frac{d}{d x}\left(\nabla_{q}\lambda_{a,m}\right)\right\rangle=0.$
    \end{enumerate}
%    \item $1$, $\left(\nabla_{q}\widetilde{\lambda}_{a,n}\right)_{n\geq 1}$
%    and $\left(\nabla_{q}\kappa_{a,n}\right)_{n\geq 1}$ generate a free family in $\lr$.
\end{enumerate}
\end{coro}

%Le résultat attendu, à savoir que les familles ci-dessous forment
%des bases de $\lr$, ne s'obtient pas aussi \og aisément \fg que
%dans le cadre sans singularités. En effet, cette propriété repose
%sur la proximité d'une base de $\lr$ au sens des asymptotiques.
%Ici, les familles de vecteurs se trouvent proches de famille de
%fonctions de Bessel dont les propriétés géométriques ne sont pas
%faciles à obtenir. Nous utiliserons pour éviter ce problème les
%opérateurs de transformation cités ci-avant.

\subsection{The spectral map}

We define the spectral map $\lambda^a\times\kappa^a:
\lr\rightarrow \R\times \ell^\infty_\R\times \ell^\infty_\R$ with
\begin{equation}\label{sl-def_lambda_x_kappa}
  \left[\lambda^a\times\kappa^a\right](q)=\left(\int_0^1 q(t)dt, \big(\widetilde{\lambda}_{a,n}(q)\big)_{n\geq 1},
   \big(n\kappa_{a,n}(q)\big)_{n\geq 1}\right    ),
\end{equation}
where $\widetilde{\lambda}_{a,n}(q)$ be given, through
\eqref{sl-integrale-In}, by
\begin{equation*}
    \lambda_{a,n}(q)=\left(n+\frac{a}{2}\right)^2\pi^2+\int_0^1 q(t)dt-a(a+1)+\widetilde{\lambda}_{a,n}(q),\quad  n\geq
    1.
\end{equation*}
Regularity for $\lambda^a\times\kappa^a$ follows as in \cite{ist}
(see also \cite{gr} when $a=1$). We give this property with the
following
\begin{thm}
    The map $\lambda^a\times\kappa^a$ is real-analytic from $\lr$ in
    $\R\times\ell^\infty_\R\times\ell^\infty_\R$. Its differential is given by
    $$d_q(\lambda^a\times\kappa^a)(v)=\Big(\left\langle 1,v\right\rangle,\big(\big\langle\nabla_{q}\widetilde{\lambda}_{a,n},v\big\rangle\big)_{n\geq 1},
    \left(\left\langle n\nabla_{q}\kappa_{a,n},v\right\rangle\right)_{n\geq 1}\Big).$$
\end{thm}
 Keep in mind that this result is \underline{not} optimal for the arrival
 space. To improve this, we use the following transformation
 operators.

\section{Transformation operators}

This idea is the key point for such an inverse spectral problem.
It was first introduced by Guillot and Ralston for $a=1$. Their
goal was to transform scalar products with the first Bessel
functions into scalar products with trigonometric functions.
Successfully used in \cite{carlson2} and \cite{cs}, Rundell and
Sacks in \cite{rusa} present theses operators stepwise for any
integer $a$: First they construct elementary operator $S_a$ which
``maps" Bessel function related to $\mathrm{H}_a(q)$ into Bessel
function related to $\mathrm{H}_{a-1}(q)$. Then, they just chain
these operators to reach Bessel function for $a=0$, in other
words, trigonometric functions. The following results extend their
results.

\begin{nota} Let $\Phi_a$ and
$\Psi_a$ be defined by
$$\Phi_a(x)=j_a(x)^2\quad\textrm{ et }\quad \Psi_a(x)=j_a(x)\eta_a(x),\quad x\in[0,1].$$
\end{nota}

Our first add, as Guillot and Ralston for $a=1$, is a
transformation property for the family $(\Psi_a)_a$ related to the
terminal velocities. A second and new improvement is some
transformation properties for the families $({\Phi_a}')_a$ and
$({\Psi_a}')_a$ related to the dual family appeared in corollary
\ref{sl-gradients_libres}. We use it further to give precision for
isospectral sets.
%%%%%%%%%%%%%%%%%%%%%%%%%%%%%%%%%%%%%%%%%%%%%%%%%%%%%%%%%%%%%%%%%%%%%%%%%

\begin{lem}%\label{lemme_reduction_indice}%
    For all $a\in\N$, $a\geq 1$ define $S_a$
    acting from $\lc$ in $\lc$ by
    \begin{equation}\label{operateur_reduction}
        S_a[f](x)=f(x)-4a\, x^{2a-1}\int_x^1
        \frac{f(t)}{t^{2a}}dt
    \end{equation}
    We have the following properties:
\begin{enumerate}[(i)]
    \item The adjoint operator for $S_a$ is \begin{equation}
    S_a^\ast[g](x)=g(x)-\frac{4a}{x^{2a}}\int_0^x t^{2a-1}
    g(t)dt.
    \end{equation}
    \item The family $\{S_a\}$ pairwise commute:%
    $$\forall (a,b)\in\N^2,\quad S_{a}S_{b}=S_{b}S_{a} .$$
    \item $S_a$ is a bounded operator in $\lc$.
    \item $S_{a}$ is a Banach isomorphism between $\lc$ and
    ${\left(x\mapsto x^{2a}\right)}^\perp$.\\
    Its inverse is given by the bounded operator on $\lc$ defined by
    $$A_{a}[g](x)=g(x)-\frac{4a}{x^{2a+1}}\int_0^x t^{2a}
    g(t)dt.$$
    \item $\Phi_a$ and $\Psi_a$ check the properties
    \begin{equation}\label{eq_reduc_indice}
        \Phi_a=-S_a^\ast[\Phi_{a-1}],\quad
        \Psi_a=-S_a^\ast[\Psi_{a-1}],
    \end{equation}
    \begin{equation}\label{eq_reduc_indice_derive}
        \Phi_a'=-A_a[\Phi_{a-1}'],\quad
        \Psi_a'=-A_a[\Psi_{a-1}'].
    \end{equation}
\end{enumerate}
\end{lem}

%%%%%%%%%%%%%%%%%%%%%%%%%%%%%%%%%%%%%%%%%%%%%%%%%%%%%%%%%%%%%%%%%%%%%%%%%

\begin{lem}\label{lemme_bessel_sinus}%
    For all $a\in\N^\ast$ let $T_a$ be defined by
    \begin{equation*}%\label{eq_lemme_bessel_sinus_def}
        T_a=(-1)^{a+1}S_a S_{a-1}\cdots S_1.
    \end{equation*}
    Then
\begin{enumerate}[(i)]
    \item $T_a$ is a bounded one-to-one operator on $\lc$
        such that :\\
         For all $q\in\lc$ and all $\lambda\in\C$,
        \begin{gather}\label{eq_lemme_bessel_sinus_1}
            \int_0^1 \big[2\Phi_a(\lambda t)-1\big]q(t)d t=\int_0^1 \cos(2\lambda t) T_a[q](t)dt,\\
        \label{eq_lemme_bessel_sinus_2}
        \int_0^1 \Psi_a(\lambda t)q(t) dt=-\frac{1}{2}\int_0^1 \sin(2\lambda t)T_a[q](t)dt.
        \end{gather}
    \item The adjoint of $T_a$ check
    \begin{equation}\label{eq_lemme_bessel_sinus_3}
        2\Phi_a(\lambda x)-1=T_a^\ast\left[\cos(2\lambda x)\right],\quad \Psi_a(\lambda x)=
        T_a^\ast\left[-\frac{1}{2}\sin(2\lambda x)\right]
    \end{equation}
    and its kernel is $$\ker(T_a^\ast)=\vect{x^2,x^4,\ldots,x^{2a}}.$$
    \item $T_a$ define a Banach isomorphism between $\lc$ and
    $\left(\ker(T_a^\ast)\right)^\perp$.\\
    Its inverse is given by the bounded operator in $\lc$ defined by
    $$B_a[f]:=(-1)^{a+1}A_a A_{a-1}\cdots A_1,$$
    moreover, we have
    \begin{equation}\label{eq_lemme_bessel_sinus_derive}
        \Phi_a'(\lambda x)=B_a\left[-\sin(2\lambda x)\right],\quad \Psi_a'(\lambda x)=
        B_a\left[-\cos(2\lambda x)\right]
    \end{equation}
\end{enumerate}
\end{lem}

%%%%%%%%%%%%%%%%%%%%%%%%%%%%%%%%%%%%%%%%%%%%%%%%%%%%%%%%%%%%%%%%%%%%%%%%%

Proof of theses results lies on properties of Bessel functions for
the calculation part and lies on Hardy inequalities for the
bounded properties and range of presented operators. For a
detailed proof, we send back to \cite{rusa}, the new facts follow
similarly.

The first use of these operators is to improve spectral data
estimations.
%%%%%%%%%%%%%%%%%%%%%%%%%%%%%%%%%%%%%%%%%%%%%%%%%%%%%%%%%%%%%%%%%%%%%%%%%
\begin{prop}
Uniformly on bounded sets in $\lr$, we have
    \begin{equation*}\label{asymptotique_spectre}
        \lambda_{a,n}(q)=\left(n+\frac{a}{2}\right)^2\pi^2+\int_0^1
        q(t)d t -a(a+1)+\ell^2(n),
    \end{equation*}
    \begin{equation*}%\label{sl-est-kappa}
        n\kappa_{a,n}(q)=\ell^2(n).
    \end{equation*}
\end{prop}

%%%%%%%%%%%%%%%%%%%%%%%%%%%%%%%%%%%%%%%%%%%%%%%%%%%%%%%%%%%%%%%%%%%%%%%%%

\begin{proof} According to \eqref{sl-integrale-In} and \eqref{sl-integrale-Jn}, we just have
to find the behavior for integral terms $I_n$ and $J_n$ defined by
$$I_n(q)=\iint_{[0,1]^2} \Big[2 j_a{\left(\omega_{a,n}(t q)x\right)^2}-1\Big]q(x)d x\,
   dt,$$
   $$J_n(q)=\iint_{[0,1]^2} \frac{j_{a}{\big(\omega_{a,n}(t q)
x\big)}\eta_a\big(\omega_{a,n}(t q)
        x\big)}{\omega_{a,n}(t q)}\,  q(x)d x\, d t.$$
We use \eqref{eq_lemme_bessel_sinus_1} and
$\eqref{eq_lemme_bessel_sinus_2}$ to get
$$I_n(q)=\iint_{[0,1]^2} \cos{\big(2\omega_{a,n}(t q)x\big)} T_a[q](x)d x\,
   dt,$$
   $$J_n(q)=-\iint_{[0,1]^2} \frac{\sin{\big(2\omega_{a,n}(t q)
x\big)}}{2\omega_{a,n}(t q)}\,  T_a[q](x)d x\, d t.$$%
From \eqref{asymptotique_spectre-bidon}, we have
$\omega_{a,n}=\left(n+\frac{a}{2}\right)\pi+\O{\frac{1}{n}}$.
Hence
$$I_n(q)=\int_0^1 \cos{\big((2n+a)\pi
x\big)} T_a[q](x)d x\,+\O{\frac{1}{n}},$$
$$J_n(q)= \frac{-1}{(2 n+a)\pi}\int_0^1\sin{\big((2n+a)\pi x\big)}\,  T_a[q](x)d x+\O{\frac{1}{n^2}}.$$
These relations lead to the result.
\end{proof}
%%%%%%%%%%%%%%%%%%%%%%%%%%%%%%%%%%%%%%%%%%%%%%%%%%%%%%%%%%%%%%%%%%%%%%%%%
Give two more estimates for the following functions
\begin{nota}
For all integer $n\geq 1$, we define
\begin{equation}\label{sl-def-Vn-Wn}
V_{a,n}(x,q):=2\left[\frac{d}{dx}\nabla_q\lambda_{a,n}\right],\quad
W_{a,n}(x,q):=-2\left[\frac{d}{dx}\nabla_q\kappa_{a,n}\right].
\end{equation}
\end{nota}
%%%%%%%%%%%%%%%%%%%%%%%%%%%%%%%%%%%%%%%%%%%%%%%%%%%%%%%%%%%%%%%%%%%%%%%%%
\begin{prop}
Uniformly the bounded sets in $[0,1]\times \lr$, we have
\begin{gather}
    \label{sl-est-Vn}V_{a,n}(x,q)=4\omega_{a,n}\left(
        {\Phi_a}'(\omega_{a,n}x)+\O{\frac{1}{n}}\right)\\
    \label{sl-est-Wn}W_{a,n}(x,q)=-2{\Psi_a}'(\omega_{a,n}x)+\O{\frac{1}{n}}.
\end{gather}
\end{prop}
%%%%%%%%%%%%%%%%%%%%%%%%%%%%%%%%%%%%%%%%%%%%%%%%%%%%%%%%%%%%%%%%%%%%%%%%%
\begin{proof}
We only proof the first one, the second follows similarly. We
have, according to the eigenvalue gradient expression,
$$V_{a,n}(x,q)=2\frac{d}{dx}\left(g_n(x,q)^2\right)=4 g_n(x,q) g_n'(x,q).$$
Estimate \eqref{sl-estimation-Rprime} combined with
\eqref{sl-preuve-est-L2-2} and \eqref{sl-def_sol_sing} leads to
$${g_n}'(x,q)=\sqrt{2} \omega_{a,n}{j_a}'(\omega_{a,n}x)+\O{1}.$$
We conclude that to \eqref{sl-est-vect-propre} which allows us to
write
$$V_{a,n}(x,q)=4\omega_{a,n} \left(2j_a(\omega_{a,n}x){j_a}'(\omega_{a,n}x)+\O{\frac{1}{n}}\right).$$
\end{proof}
%%%%%%%%%%%%%%%%%%%%%%%%%%%%%%%%%%%%%%%%%%%%%%%%%%%%%%%%%%%%%%%%%%%%%%%%%

        %%%%%%%%%%%%%%%%%%%%%%%%%%%%%%%%%%%%%%%%%%%%%%%%%%%%%%%%%%%%%%%%%%%%
%%%%%%%%%%%%%%%%%%%%%%%%%%%%%%%%%%%%%%%%%%%%%%%%%%%%%%%%%%%%%%%%%%%%
\section{The inverse spectral problem}
%%%%%%%%%%%%%%%%%%%%%%%%%%%%%%%%%%%%%%%%%%%%%%%%%%%%%%%%%%%%%%%%%%%%

Now, the map $\lambda^a\times\kappa^a$ defined by
\eqref{sl-def_lambda_x_kappa} is a real analytic map from $\lr$
into $\R\times\ell^2_\R\times\ell^2_\R$. Give the main result of
this paper:
%%%%%%%%%%%%%%%%%%%%%%%%%%%%%%%%%%%%%%%%%%%%%%%%%%%%%%%%%%%%%%%%%%%%
\begin{thm}\label{th_diff_inversible}
    For all integer $a\geq 1$, $d_q(\lambda^a\times\kappa^a)$ is a Banach isomorphism from $\lr$ onto
    $\R\times\ell^2_\R\times\ell^2_\R$.
\end{thm}

From the results of Zhornitskaya, Zerov \cite{zs} and Carlson
\cite{carlson}, for all $a\geq -\frac{1}{2}$, the map
$\lambda^a\times\kappa^a$ is one-to-one on $\lr$. We immediately
deduce the following property

%%%%%%%%%%%%%%%%%%%%%%%%%%%%%%%%%%%%%%%%%%%%%%%%%%%%%%%%%%%%%%%%%%%%

\begin{coro} For all integer $a\geq 1$, the map $\lambda^a\times\kappa^a$ is a global
real analytic coordinate system on $\lr$.
\end{coro}

%%%%%%%%%%%%%%%%%%%%%%%%%%%%%%%%%%%%%%%%%%%%%%%%%%%%%%%%%%%%%%%%%%%%
\begin{proof}[Proof of the theorem] The main arguments of this proof are extensions of those in \cite{gr}.
From orthogonality relations in corollary
\ref{sl-gradients_libres} the family
$$\left\{1\right\}\cup\left\{\nabla_{q}\widetilde{\lambda}_{a,n}\right\}_{n\geq 1}\cup \left\{\nabla_q\kappa_{a,n}\right\}_{n\geq 1}$$
is free in $\lr$. Let $r_n$ and $s_n$ be defined by
\begin{gather}
\label{sl-def-rn}  r_n(x)=\nabla_{q}\widetilde{\lambda}_{a,n}(x)-\left(2\Phi_a(\omega_{a,n} x)-1\right),\\
\label{sl-def-sn}
s_n(x)=\nabla_{q}\kappa_{a,n}(x)-\frac{1}{\omega_{a,n}}\Psi_a(\omega_{a,n}
x).
\end{gather}
Using lemma \ref{lemme_bessel_sinus}, we have, for all $v\in\lr$,
\begin{gather}
    \label{eq_grad_pert_1}\left\langle\nabla_{q}\widetilde{\lambda}_{a,n},v\right\rangle=\int_0^1
    \left(\cos{(2\omega_{a,n} t)}+R_n(t)\right) T_a[v](t)dt,\\
    \label{eq_grad_pert_2}\left\langle\nabla_{q}\kappa_{a,n},v\right\rangle=\frac{-1}{2\omega_{a,n}}\int_0^1
    \left(\sin{(2\omega_{a,n} t)}+S_n(t)\right)T_a[v](t)dt,
\end{gather}
with
\begin{equation}\label{sl-def-Rn-Sn}
    R_n=B_a^\ast[r_n]\textrm{ et }
\frac{-S_n}{2\omega_{a,n}}=B_a^\ast[s_n].
\end{equation}
 Moreover, recall that
$T_a[1]=-1$, particularly
\begin{equation*}%\label{eq_grad_pert_3}
\left\langle1,v\right\rangle=\int_0^1
v(t)dt= \int_0^1 -T_a[v](t)dt.
\end{equation*}
We denote $F$ the operator defined by
\begin{multline*}
   F(w)=\bigg(\left\langle-1,w\right\rangle,\left\{
        \left\langle \cos{(2\omega_{a,n} t)}+R_n(t),w\right\rangle
                \right\}_{n\geq 1},\\\left\{
        \left\langle\frac{-n}{2\omega_{a,n}}\left(\sin{(2\omega_{a,n} t)}+S_n(t)\right),w\right\rangle
                \right\}_{n\geq 1}\bigg),
\end{multline*}
such that $d_{q}(\lambda^a\times\mu^a)(v)=F\circ T_a[v]$. From
lemma \ref{lemme_bessel_sinus}, $T_a$ is invertible from $\lr$
onto
$\left(\displaystyle\vect{x^2,x^4,\ldots,x^{2a}}\right)^\perp$.
Thus, we have to show that $F$ is invertible from
$\left(\displaystyle\vect{x^2,x^4,\ldots,x^{2a}}\right)^\perp$
onto $\R\times\ell^2_\R\times\ell^2_\R$. Following \cite{gr}, we
have to show the invertibility of $\mathbf{F}$ the operator from
$\lr$ into $\R \times\ell^2_\R\times\ell^2_\R$, sending functions
in $\lr$ on their Fourier coefficients (or, more simply, their
scalar products against each element of the considered family)
with respect to the family
\begin{multline}\label{sl-def-famille-libre}
    \mathcal{F}=\left\{1\right\}\cup\left\{t^{2j}\right\}_{j\in\ing 1,a\ind}\cup\left\{
            \cos{(2\omega_{a,n} t)}+R_n(t)\right\}_{n\geq 1}\\
            \cup
            \left\{\frac{-n}{2\omega_{a,n}}\left(\sin{(2\omega_{a,n} t)}+S_n(t)\right)\right\}_{n\geq
            1}.
\end{multline}
This last statement will come from \cite{ist}: Appendix D, theorem
3:

\begin{lem} Let $\mathcal{F}=\{f_n\}_{n\in\N}$ be a sequence
in an Hilbert space $H$, with the two following properties
\begin{enumerate}[(a)]
    \item there exists an orthogonal basis $\mathcal{E}=\{e_n\}_{n\in\N}$ in $H$
    for which $$\sum \n{f_n-e_n}{2}^2<\infty.$$
    \item the $f_n$ are linearly independent.
\end{enumerate}
Then $\left\{f_n\right\}_{n\in\N}$ is a basis of $H$ and the map
$\mathbf{F}:x\mapsto\{(f_n,x)\}_{n\in\N}$ is a linear isomorphism
between $H$ and $\ell^2_\R$.
\end{lem}
%%%%%%%%%%%%%%%%%%%%%%%%%%%%%%%%%%%%%%%%%%%%%%%%%%%%%%%%%%%%%%%%%%%%%
\noindent Considering the orthonormal basis
    $$\mathcal{E}=\left\{\sqrt{2}\cos{\pi x},\sqrt{2}\sin{\pi
    x},\ldots,\sqrt{2}\cos{(2n+1)\pi x},\sqrt{2}\sin{(2n+1)\pi
    x},\ldots\right\},$$
    when $a=2\mathfrak{a}+1$ with $\mathfrak{a}\in\N$ (as in \cite{gr} for $a=1$),
    and
        $$\mathcal{E}=\left\{1,\sqrt{2}\cos{2\pi x},\sqrt{2}\sin{2\pi
    x},\ldots,\sqrt{2}\cos{2n\pi x},\sqrt{2}\sin{2n\pi x},\ldots\right\},$$
    when $a=2\mathfrak{a}$ with $\mathfrak{a}\in\N$ (as in \cite{ist} for $a=0$), estimates \eqref{gradient_spectre_estimation}
    and \eqref{gradient_kappa_estimation}
    used with \eqref{sl-def-rn} and \eqref{sl-def-sn} show that
    $$\n{r_n}{\lr}=\O{\frac{1}{n}},\quad \n{s_n}{\lr}=\O{\frac{1}{n^2}},$$
    then the boundedness of $B_a^\ast$ and relations given in
    \eqref{sl-def-Rn-Sn} give
    $$\n{R_n}{\lr}=\O{\frac{1}{n}},\quad \n{S_n}{\lr}=\O{\frac{1}{n}},$$
    which implies condition $(a)$ after normalization. Lemma
    \ref{lemme_famille_libre} will give condition $(b)$, thus the proof will be complete.
\end{proof}
%%%%%%%%%%%%%%%%%%%%%%%%%%%%%%%%%%%%%%%%%%%%%%%%%%%%%%%%%%%%%%%%%%%%

Now, give the lemma implying the theorem. Ingredients of the proof
are similar to \cite{gr}, we improve the method in a systematic
way (for each $a$).
\begin{lem} \label{lemme_famille_libre} $\mathcal{F}$ defined by \eqref{sl-def-famille-libre} is a free family in $\lr$.
\end{lem}

%%%%%%%%%%%%%%%%%%%%%%%%%%%%%%%%%%%%%%%%%%%%%%%%%%%%%%%%%%%%%%%%%%%%

\begin{proof}  Rewrite relations \eqref{eq_grad_pert_1} and
\eqref{eq_grad_pert_2} as
$$T_a^\ast\left(\cos{(2\omega_{a,n}
t)}+R_n(t)\right)=\nabla_{q}\widetilde{\lambda}_{a,n}$$%
and
$$ T_a^\ast\left(\frac{-1}{2\omega_{a,n}}\left(\sin{(2\omega_{a,n}
t)}+S_n(t)\right)\right)=\nabla_{q}\kappa_{a,n}.$$%
Since $T_a^\ast$ is bounded and the famille
$\left\{1\right\}\cup\left\{\nabla_q\widetilde{\lambda}_{a,n}\right\}_{n\geq
1}\cup\left\{\nabla_q\kappa_{a,n}\right\}_{n\geq 1}$ is free, we
deduce the linear independency for the following family
$$\left\{1\right\}\cup\left\{
            \cos{(2\omega_{a,n} t)}+R_n\right\}_{n\geq 1}\\
                \cup\left\{\frac{-1}{2\omega_{a,n}}\left(\sin{(2\omega_{a,n} t)}+S_n\right)\right\}_{n\geq
                1}.$$
Let $k\in \ing 1,a\ind$, denote $W_k$ the function
$W_{k}(t)=t^{2k}$. Let us show that $W_k$ is not in the closure of
$\operatorname{Vect}\left(\mathcal{F}\setminus\left\{W_k\right\}\right)$.
(We might show iteratively that
$W_{k}\notin\overline{\operatorname{Vect}\left(\mathcal{F}\setminus\left\{W_j,j\in\ing
k,a\ind\right\}\right)}$, but it is not necessary by taking $\alpha_m^{(j)}=0$ for any $m\in\ing k,a\ind$ in the next expression.)\\
Suppose the contrary: there exists for $j\in\N$ a sequence of
vector
\begin{multline*}
    W_k^{(j)}(t)=\alpha_0^{(j)}+\sum_{\substack{m\in\ing 1,a\ind\\m\neq k}} \alpha_m^{(j)}\,W_m(t)+\sum_{n\in\ing 1,N_j\ind}a_n^{(j)}\left(
            \cos{(2\omega_{a,n} t)}+R_n(t)\right) \\
    +\sum_{n\in\ing 1,N_j\ind}b_n^{(j)}\frac{-1}{2\omega_{a,n}}\left(\sin{(2\omega_{a,n} t)}+S_n(t)\right),
\end{multline*}
with $N_j<\infty$, $ \alpha_m^{(j)},a_n^{(j)},b_n^{(j)}\in\R$ such
that
$$W_k^{(j)}\underset{j\rightarrow\infty}\longrightarrow W_k \textrm{ in }
\lr.$$%
Since $T_a^\ast(W_m)=0$ for $m=1,\dots,a$, the sequence
$$w^{(j)}:=T_a^\ast(W_k^{(j)})=-\alpha_0^{(j)}+\sum_{n\in\ing 1,N_j\ind}a_n^{(j)}\nabla_{q}\widetilde{\lambda}_{a,n}
+b_n^{(j)}\nabla_{q}\kappa_{a,n}$$ tends to $0$ in $\lr$ when
$j\rightarrow\infty$. Thus, corollary \ref{sl-gradients_libres}
gives
\begin{gather}
    \label{aoj_zero}\alpha_0^{(j)}=\int_0^1 w^{(j)}(t) d t\,\underset{j\rightarrow\infty}{\longrightarrow}0,\\
    \label{anj_zero}a_n^{(j)}=-2\int_0^1 w^{(j)}(t)\frac{d}{d x}\left(\nabla_q\kappa_{a,n}\right) d t\,\underset{j\rightarrow\infty}{\longrightarrow}0,\\
    \label{bnj_zero}b_n^{(j)}=-2\int_0^1 w^{(j)}(t)\frac{d}{d
    x}\left(\nabla_q \lambda_{a,n}\right) d t\,\underset{j\rightarrow\infty}{\longrightarrow}0.
\end{gather}

%%%%%%%%%%%%%%%%%%%%%%%%%%%%%%%%%%%%%%%%%%%%%%%%%%%%%%%%%%%%%%%%%%%%
%%%% INTRODUCTION DE LA FONCTION TEST
%%%%%%%%%%%%%%%%%%%%%%%%%%%%%%%%%%%%%%%%%%%%%%%%%%%%%%%%%%%%%%%%%%%%
Now consider $\omega\in \classe{\infty}_0([0,1],\R)$, supported in
$[\delta,1]$ with $\delta>0$, such that
$$\left\langle\omega, W_m\right\rangle=\delta_{k,m},\quad m\in\ing 1,a\ind$$
and $$\left\langle B_a[\omega], 1\right\rangle=0
\quad\textrm{i.e.}\quad \left\langle\omega,
1\right\rangle=c(k,a)\left\langle\omega,
W_k\right\rangle.%  Note de bas de page
\footnote{$c(k,a)$
 represents the $W_k$ component of
$B_a^\ast[1]$. It can by computed by induction on $a$:
$B_a^\ast[1]=-1+\sum_{m=1}^a c(m,a)W_m.$
Particularly $c(1,a)=a(a+1)$ (idem \cite{gr} when $a=1$).}%   Fin de la note
$$
 Smoothness and support of $\omega$ imply that
$$\int_0^1 \omega(t) \cos{(2\omega_{a,n} t)}d t,\int_0^1 \omega(t)\sin{(2\omega_{a,n} t)}d t=\O{\frac{1}{n^N}},
\quad \forall N\in\N$$ and that $B_a[\omega]$ is
$\classe{\infty}([0,1],\R)$ supported in $[\delta,1]$.\\
Now, plug estimation \eqref{sl-estimation-R-L2} in the integral
expression \eqref{sl-eq-int-phi}, then use
\eqref{sl-preuve-est-L2-2}, controls \eqref{annexe-bessel_est_ja}
and \eqref{annexe-resolvante-sl_est_tx} to obtain the uniform
estimate on $[0,1]$
\begin{multline}\label{sl-exp-rn}
  r_n(x)= \frac{1}{\n{\varphi_n}{2}^2}\left(2 u(x,\lambda_{a,n})\int_0^x \G(x,t,\lambda_{a,n})q(t)u(t,\lambda_{a,n})dt\right)\\
   +\left(\frac{1}{\n{\varphi_n}{2}^2}-\frac{2
    {\left(\lambda_{a,n}\right)}^{a+1}}{((2a+1)!!)^2}\right)u(x,\lambda_{a,n})^2+\O{\frac{1}{n^2}}.
\end{multline}
Thanks to lemma \ref{lemme-sommable-1}, this implies that the
following family is summable
$$\left\{\left\langle\omega,\big[\cos{(2\omega_{a,n} t)}+R_n(t)\big]\right\rangle\right\}_{n\geq 1}.$$

%%%%%%%%%%%%%%%%%%%%%%%%%%%%%%%%%%%%%%%%%%%%%%%%%%%%%%%%%%%%%%%%%%%%%%%%
Now turn to $s_n$. From relations \eqref{sl-def-sn} and
\eqref{gradient_kappa}, we get
\begin{eqnarray*}
  s_n(x) &=& -a_n(x,q)+\nabla_q\lambda_{a,n}(x)\int_0^1 a_n(t,q)dt-\frac{1}{\omega_{a,n}}\Psi_a(\omega_{a,n}x), \\
         &=& -a_n(x,q)+\left(2\Phi_a(\omega_{a,n}x)+r_n(x)\right)\int_0^1 a_n(t,q)dt-\frac{1}{\omega_{a,n}}\Psi_a(\omega_{a,n}x), \\
         &=& -a_n(x,q)+2\Phi_a(\omega_{a,n}x)\int_0^1 a_n(t,q)dt \\
         & &\phantom{ -a_n(x,q)}+r_n(x)\int_0^1
         a_n(t,q)dt-\frac{1}{\omega_{a,n}}\Psi_a(\omega_{a,n}x).
\end{eqnarray*}
Respectfully insert \eqref{sl-estimation-R-L2} and
\eqref{sl-estimation-S-L2} in integral expressions
\eqref{sl-eq-int-phi} and \eqref{sl-eq-int-psitilde}; then use
\eqref{annexe-bessel_est_ja}, \eqref{annexe-bessel_est_na},
\eqref{annexe-resolvante-sl_est_tx},
\eqref{annexe-resolvante-sl_est_xt} and wronskian estimate
\eqref{sl-def_sol_sing} to obtain the uniform estimate on $[0,1]$
\begin{equation}\label{def-exp-sn}
\begin{aligned}&\left.
\begin{aligned}
 s_n(x)= & -v(x,\lambda_{a,n})\int_0^x \G(x,t,\lambda_{a,n})q(t)u(t,\lambda_{a,n})dt\\
    &+u(x,\lambda_{a,n})\int_x^1 \G(x,t,\lambda_{a,n})q(t)v(t,\lambda_{a,n})dt\\
    &+2\Phi_a(\omega_{a,n}x)\int_0^1
   a_n(t,q)dt\\
\end{aligned}\right\}\tilde{s}_n(x)\\
   &\qquad\qquad+r_n(x)\int_0^1
   a_n(t,q)dt\\
   &\qquad\qquad+\frac{1}{\omega_{a,n}}\Psi_a(\omega_{a,n}x)\left(W^{-1}-1\right)+\O{\frac{1}{n^3}}.
\end{aligned}
\end{equation}
With the help of lemma \ref{lemme-sommable-2}, we deduce the
summability of the family
$$\left\{n\int_0^1\omega(t)\frac{-1}{2\omega_{a,n}}\left(\sin{(2\omega_{a,n} t)}+S_n(t)\right)dt\right\}_{n\geq 1}.$$
We finish the proof writing
\begin{multline*} \left\langle\omega,
W_k^{(j)}\right\rangle=\alpha_0^{(j)}\left\langle\omega,
1\right\rangle+\sum_{n\in\ing
1,N_j\ind}a_n^{(j)}\left\langle\omega,\left(
    \cos{(2\omega_{a,n} t)}+R_n(t)\right)\right\rangle\\
    +\sum_{n\in\ing 1,N_j\ind}b_n^{(j)}\left\langle\omega,\frac{-1}{2\omega_{a,n}}\left(\sin{(2\omega_{a,n}
    t)}+S_n(t)\right)\right\rangle.
\end{multline*}
Indeed, the following estimates, deduced from
\eqref{sl-est-Vn}-\eqref{sl-est-Wn} and
\eqref{aoj_zero}-\eqref{bnj_zero},
$$\big|a_n^{(j)}\big|\leq C,\quad\big|b_n^{(j)}\big|\leq C n,\quad\forall (n,j)\in\N^2,$$ imply that
$$\left\langle\omega,
W_k^{(j)}\right\rangle\underset{j\rightarrow\infty}\longrightarrow
0,$$ which is  contradictory with the choice of $\omega$.
\end{proof}

%%%%%%%%%%%%%%%%%%%%%%%%%%%%%%%%%%%%%%%%%%%%%%%%%%%%%%%%%%%%%%%%%%%%%%%%%%%%%%%%%%

\begin{lem}\label{lemme-sommable-1}
$$\left\{\left\langle\omega,R_n\right\rangle\right\}_{n\geq 1}\in\ell^1_\R.$$
\end{lem}

%%%%%%%%%%%%%%%%%%%%%%%%%%%%%%%%%%%%%%%%%%%%%%%%%%%%%%%%%%%%%%%%%%%%%%%%%%%%%%%%%%

\begin{proof}
Let $r_{n,1}$ and $r_{n,2}$ be the first and second terms
in\eqref{sl-exp-rn}. Recall that
$$\langle\omega,R_n\rangle=\langle\omega,B_a^\ast[r_n]\rangle=\langle B_a[\omega],r_n\rangle,$$
thus, we just have to show that  $\{\langle
B_a[\omega],r_{n,j}\rangle\}\in\ell^1_\R$, $j=1,2$. First consider
$r_{n,2}$:
\begin{eqnarray*}
  \langle
B_a[\omega],r_{n,2}\rangle &=& \left\langle
B_a[\omega],\left(\frac{1}{\n{\varphi_n}{2}^2}-\frac{2
    {\left(\lambda_{a,n}\right)}^{a+1}}{((2a+1)!!)^2}\right)u(x,\lambda_{a,n})^2\right\rangle, \\
   &=&\left\langle
B_a[\omega],\left(\frac{((2a+1)!!)^2}{2
    {\left(\lambda_{a,n}\right)}^{a+1}\n{\varphi_n}{2}^2}-1\right)2
    j_a(\omega_{a,n}x)^2\right\rangle.
\end{eqnarray*}
With \eqref{sl-preuve-est-L2-2}, we have
$$  \langle
B_a[\omega],r_{n,2}\rangle =
\O{\frac{1}{\omega_{a,n}}}\left\langle B_a[\omega],2
j_a(\omega_{a,n}x)^2\right\rangle,
$$
then, hypothesis upon $B_a[\omega]$ gives
$$  \langle
B_a[\omega],r_{n,2}\rangle =
\O{\frac{1}{\omega_{a,n}}}\left\langle B_a[\omega],2
j_a(\omega_{a,n}x)^2-1\right\rangle.
$$
Lemma \ref{lemme_bessel_sinus} leads to the result.\\
Now consider $r_{n,1}$ and write as before
\begin{eqnarray*}
r_{n,1}(x) &=&  \frac{1}{\n{\varphi_n}{2}^2}
\left(2 u(x,\lambda_{a,n})\int_0^x \G(x,t,\lambda_{a,n})q(t)u(t,\lambda_{a,n})dt\right) \\
   &=&\frac{2((2a+1)!!)^2}{(\omega_{a,n})^{2a+2}\n{\varphi_n}{2}^2}
\, j_a(\omega_{a,n}x)\int_0^x
\G(x,t,\lambda_{a,n})q(t)j_a(\omega_{a,n}t)dt.
\end{eqnarray*}
With \eqref{annexe-bessel_est_ja} and
\eqref{annexe-resolvante-sl_est_tx} we get
\begin{equation*}%\label{eq-bidon-11}
    \int_0^x
\G(x,t,\lambda_{a,n})q(t)j_a(\omega_{a,n}t)dt=\O{\frac{1}{\omega_{a,n}}},
\end{equation*}
and, thanks to \eqref{sl-preuve-est-L2-2}, we have
$$r_{n,1}(x)=4 j_a(\omega_{a,n}x)\int_0^x
\G(x,t,\lambda_{a,n})q(t)j_a(\omega_{a,n}t)dt+\O{\frac{1}{n^2}}.$$
From the expression of $\G(x,t,\lambda)$, we can write
\begin{align*}
    \langle B_a[\omega],& \,r_{n,1}\rangle=\\
    &\frac{4}{\omega_{a,n}} \int_0^1 B_a[\omega](x)j_a(\omega_{a,n}x)^2\int_0^x
q(t)j_a(\omega_{a,n}t)\eta_a(\omega_{a,n}t)d t d
x\\
&-\frac{4}{\omega_{a,n}} \int_0^1
B_a[\omega](x)j_a(\omega_{a,n}x)\eta_a(\omega_{a,n}x)\int_0^x
q(t)j_a(\omega_{a,n}t)^2d t d x\\
 &+\O{\frac{1}{n^2}}.
\end{align*}
But, $B_a[\omega]$ is supported in $[\delta,1]$, $\delta>0$.
Relations \eqref{annexe-bessel-ja-sinus} and
\eqref{annexe-bessel-na-sinus} thus give
\begin{align*}
    \langle B_a[\omega],& \,r_{n,1}\rangle=\\
    &\frac{2}{\omega_{a,n}} \int_0^1 B_a[\omega](x)(1-\cos(2\omega_{a,n}x-a\pi))\int_0^x
q(t)j_a(\omega_{a,n}t)\eta_a(\omega_{a,n}t)d t d
x\\
&-\frac{2}{\omega_{a,n}} \int_0^1
B_a[\omega](x)\sin(2\omega_{a,n}x-a\pi)\int_0^x
q(t)j_a(\omega_{a,n}t)^2d t d x\\
 &+\O{\frac{1}{n^2}}.
\end{align*}
Integrating by parts in all terms having factors of
$\cos(\omega_{a,n}x)$ or $\sin(\omega_{a,n}x)$, we get
\begin{align*}
    \langle B_a[\omega], \,r_{n,1}\rangle=&\frac{2}{\omega_{a,n}} \int_0^1 B_a[\omega](x)\int_0^x
q(t)j_a(\omega_{a,n}t)\eta_a(\omega_{a,n}t)d t d
x\\
 &+\O{\frac{1}{n^2}}.
\end{align*}
Finally, interchanging the order of integration and using
properties of the transformation operator, we obtain the
concluding relation
$$
    \langle B_a[\omega], \,r_{n,1}\rangle=\frac{-1}{\omega_{a,n}} \int_0^1
    \sin(2\omega_{a,n}t) T_a\left[t\mapsto
    q(t)\int_t^1 B_a[\omega](x) d x\right] d t+\O{\frac{1}{n^2}}.
$$
\end{proof}

%%%%%%%%%%%%%%%%%%%%%%%%%%%%%%%%%%%%%%%%%%%%%%%%%%%%%%%%%%%%%%%%%%%%%%%%%%%%%%%%%%

\begin{lem}\label{lemme-sommable-2}
$$\left\{\left\langle\omega,S_n\right\rangle\right\}_{n\geq 1}\in \ell^1_\R.$$
\end{lem}

%%%%%%%%%%%%%%%%%%%%%%%%%%%%%%%%%%%%%%%%%%%%%%%%%%%%%%%%%%%%%%%%%%%%%%%%%%%%%%%%%%

\begin{proof}
Recall
$$\langle\omega, S_n\rangle=2\omega_{a,n}\langle\omega,B_a^\ast[s_n]\rangle=2\omega_{a,n}\langle B_a[\omega],s_n\rangle.$$
The proof is similar to lemma \ref{lemme-sommable-1} for
$\tilde{s}_n$ (see \eqref{def-exp-sn}): we change in the scalar
product $\langle B_a[\omega],\tilde{s}_n\rangle$,
$2\Phi_a(\lambda_{a,n}x)$ by $2\Phi_a(\lambda_{a,n}x)-1$ and use
the transformation operator; for the two following terms in
$\tilde{s}_n$, as previously, we have to use the support of
$B_a[\omega]$, asymptotics deduced from
\eqref{annexe-bessel-ja-sinus} and \eqref{annexe-bessel-na-sinus},
integrate by parts terms with $\cos{(2\omega_{a,n}x)}$ or
$\sin{(2\omega_{a,n}x)}$ and finally invert integration order to
use transformation operators.

The term following $\tilde{s}_n $ in \eqref{def-exp-sn} is
controlled thanks to $r_n$; denote $\hat{s}_n$ the remaining term.
But, we have
\begin{eqnarray*}
  \left\langle B_a[\omega],\hat{s}_n\right\rangle &=& \O{\frac{1}{n^2}}\int_0^1 B_a[\omega](t)\ja{a}{_{a,n} t}\na{a}{_{a,n} t}dt \\
   &=&\O{\frac{1}{n^2}}\int_\delta^1 B_a[\omega](t)\sin{(2\omega_{a,n}
   t)}dt+\O{\frac{1}{n^3}}.
\end{eqnarray*}
Thus $\Big\{n\left\langle
B_a[\omega],\hat{s}_n\right\rangle\negthinspace\Big\}_{n\geq 1}$
is in $\ell^1_\R$ and the proof is completed.
\end{proof}

%%%%%%%%%%%%%%%%%%%%%%%%%%%%%%%%%%%%%%%%%%%%%%%%%%%%%%%%%%%%%%%%%%%%%%%%%%%%%%%%%%

Now comes the corollary (see \eqref{sl-def-Vn-Wn} for the
definition of $V_n$ and $W_n$).

%%%%%%%%%%%%%%%%%%%%%%%%%%%%%%%%%%%%%%%%%%%%%%%%%%%%%%%%%%%%%%%%%%%%%%%%%%%%%%%%%%

\begin{coro}\label{sl-exp-inverse-diff} $\lambda^a\times\kappa^a$ is a local real analytic isomorphism on $\lr$.
Moreover, the inverse of $d_q(\lambda^a\times\kappa^a)$ is the
bounded linear map from $\R\times\ell^2_\R\times\ell^2_\R$ onto
$\lr$ given by
$$(d_q(\lambda^a\times\kappa^a))^{-1}(\eta_0,\eta,\xi)=\eta_0+\sum_{n\geq 1}\eta_n W_{a,n}+\sum_{n\geq 1}\frac{\xi_n}{n} V_{a,n}.$$
\end{coro}

%%%%%%%%%%%%%%%%%%%%%%%%%%%%%%%%%%%%%%%%%%%%%%%%%%%%%%%%%%%%%%%%%%%%%%%%%%%%%%%%%%

\begin{proof}The first part of the corollary comes from the theorem \ref{th_diff_inversible} and from
the inverse function theorem. Now consider
$(\eta_0,\eta,\xi)\in\R\times\ell^2_\R\times\ell^2_\R$ and define%
$$u=\eta_0+\sum_{n\geq 1}\eta_n W_{a,n}+\sum_{n\geq 1}\frac{\xi_n}{n} V_{a,n}.$$
Since $B_a$ is bounded, estimations \eqref{sl-est-Vn} and
\eqref{sl-est-Wn}, with relations
\eqref{eq_lemme_bessel_sinus_derive} give
$$\frac{1}{n}V_{a,n}(x,q)=B_a\left[\frac{4\omega_{a,n}}{n} \sin(2\omega_{a,n}x)+\O{\frac{1}{n}}\right]$$
and
$$W_{a,n}(x,q)=B_a\left[-2 \cos(2\omega_{a,n}x)+\O{\frac{1}{n}}\right].$$

Definition of $\xi$ and $\eta$, estimation of eigenvalues and
boundedness of $B_a$ imply the convergence in $\lr$ for the series
defining $u$.

From corollary \ref{sl-gradients_libres}, we have
$$\left\langle 1,u\right\rangle=\eta_0,$$ and for all integer $n\geq 1$%
$$\left\langle \nabla_{q}\widetilde{\lambda}_{a,n},u\right\rangle=\eta_n, \quad \left\langle n\nabla_{q}\kappa_{a,n},u\right\rangle=\xi_n.$$
Thus, we have
$d_{q}\left(\lambda^a\times\kappa^a\right)(u)=(\eta_0,\eta,\xi)$,
which proofs the corollary.
\end{proof}
%%%%%%%%%%%%%%%%%%%%%%%%%%%%%%%%%%%%%%%%%%%%%%%%%%%%%%%%%%%%%%%%%%%%%%%%%%%%%%%%%%
We finish with the description of isospectral sets. For
$q_0\in\lr$, we define the set of potentials with same Dirichlet
spectrum as $q_0$, called isospectral set, by
$\operatorname{Iso}(q_0,a)=\left\{q\in\lr:
\lambda^a(q)=\lambda^a(q_0)\right\}.$ The new fact of the
following result is to explicit tangent and normal spaces.
%%%%%%%%%%%%%%%%%%%%%%%%%%%%%%%%%%%%%%%%%%%%%%%%%%%%%%%%%%%%%%%%%%%%%%%%%%%%%%%%%%
\begin{thm} Let $q_0\in\lr$, then
\begin{enumerate}[(a)]
    \item $\operatorname{Iso}(q_0,a)$ is a real-analytic manifold of $\lr$ of infinite dimension and codimension, lying in the hyperplane
    of all functions with mean $\int_0^1 q_0(t)dt$.
    \item At every point $q$ in $\operatorname{Iso}(q_0,a)$, the tangent
    space is
    $$T_q\operatorname{Iso}(q_0,a)=\left\{\sum_{n\geq 1}\frac{\xi_n}{n}V_{a,n}: \xi\in\ell^2_\R\right\}$$
    and the normal space is
    $$N_q\operatorname{Iso}(q_0,a)=\left\{\eta_0+\sum_{n\geq 1}\eta_n ({g_n}^2-1): (\eta_0,\eta)\in\R\times\ell^2_\R\right\}.$$
\end{enumerate}
\end{thm}
%%%%%%%%%%%%%%%%%%%%%%%%%%%%%%%%%%%%%%%%%%%%%%%%%%%%%%%%%%%%%%%%%%%%%%%%%%%%%%%%%%
\begin{proof}
It is straightforward from \cite{ist}.
\begin{enumerate}[(a)]
    \item % Let $q_0\in\lr$ and
%    $q\in\operatorname{Iso}(q_0,a)$. From corollary
%    \ref{sl-exp-inverse-diff}, the $\lambda^a\times\kappa^a: \lr\rightarrow\R\times\ell^2_\R\times\ell^2_\R$
%    is a real-analytic diffeomorphism locally in a neighborhood of $q$, that is there exists an open set
%    $\mathcal{V}\subset\lr$ containing $q$ and an open set
%    $\mathcal{W}\subset\R\times\ell^2_\R\times\ell^2_\R$ such that
%    $\lambda^a\times\kappa^a: \mathcal{V}\rightarrow\mathcal{W}$
%    is a real-analytique isomorphisme. Then, we have the identification
%    $$\left(\lambda^a\times\kappa^a\right)(\mathcal{V}\cap\operatorname{Iso}(q_0))=
%    \mathcal{W}\cap\left(\left\{\lambda^a(q_0)\right\}\times \ell^2_\R\right),$$
%    which shows that $\operatorname{Iso}(q_0)$ is a real-analytic
%    manifold in $\lr$.
    The first part of the assertion comes from \cite{carlson2}
    Theorem 1.3, the second is direct from
    \eqref{asymptotique_spectre}.
    \item Since $$T_{q}\operatorname{Iso}(q_0,a)=
    \left(d_{q}\left(\lambda^a\times\kappa^a\right)\right)^{-1}\big(\big\{0_{\R\times\ell^2_\R}\big\}\times\ell^2_\R\big),$$
    corollary \ref{sl-exp-inverse-diff} gives expression of the tangent space.
    According to corollary \ref{sl-gradients_libres}, $\left\{1, {g_n}^2-1 :n\geq 1\right\}$ is free, orthogonal
    to the free family $(V_{a,n})_{n\in\Z}$. Thus, we have
    $$\left\{\eta_0+\sum_{n\geq 1}\eta_n (g_n^2-1): (\eta_0,\eta)\in\R\times\ell^2_\R\right\}\subset N_q\operatorname{Iso}(q_0).$$
    Moreover, any vector orthogonal to $\left\{1, {g_n}^2-1 :n\geq
    1\right\}$ is, with regards to the Fréchet derivative of $\lambda^a$, in the null space of $d_q\lambda^a$. Thus, we get the other inclusion and
    then the proof.
\end{enumerate}
\end{proof}
%%%%%%%%%%%%%%%%%%%%%%%%%%%%%%%%%%%%%%%%%%%%%%%%%%%%%%%%%%%%%%%%%%%%%%%%%%%%%%%%%%
To finish, we recall that the characterization of the spectra of
each operator $\mathrm{H}_a(q)$ was obtained by Carlson
(\cite{carlson2} Theorem 1.1).

%   ANNEXES
        %   Annexe A: Fonctions de Bessel
\section{Bessel functions} \label{annexe-bessel}

Spherical Bessel functions $j_a$ and $\eta_a$ are defined through
\begin{equation}\label{annexe-bessel-def}
    j_a(z)=\sqrt{\frac{\pi z}{2}}J_{a+1/2}(z),\quad
    \eta_a(z)=(-1)^a \sqrt{\frac{\pi z}{2}}J_{-a-1/2}(z),
\end{equation}
where $J_\nu$ is the first kind Bessel Bessel function of order
$\nu$ (see \cite{emot} for precisions).

From \cite{emot} formulas $(1)$ and $(2)$ section $7.11$ p.78,
they behave like
    \begin{equation}\label{annexe-bessel-ja-sinus}
    j_a(z)=\sin{\left(z-\frac{a\pi}{2}\right)}+\O{\frac{e^{|\im z|}}{|z|}}, \quad |z|\rightarrow\infty,
    \end{equation}
    \begin{equation}\label{annexe-bessel-na-sinus}
    \eta_{a}(z)=\cos{\left(z-\frac{a\pi}{2}\right)}+\O{\frac{e^{|\im
    z|}}{|z|}}, \quad |z|\rightarrow\infty.
    \end{equation}
The following estimates can be found in \cite{carlson}
\begin{itemize}
    \item Uniform estimates on $\C$:
    \begin{eqnarray}
      \label{annexe-bessel_est_ja}\left|j_{a}(z)\right| &\leq& C e^{|\im z|}\left(\frac{|z|}{1+|z|}\right)^{a+1}, \\
      \label{annexe-bessel_est_na}\left|\eta_{a}(z)\right| &\leq& C e^{|\im z|}\left(\frac{1+|z|}{|z|}\right)^{a}.
    \end{eqnarray}
    \item Estimations for the Green function $G(x,t,\lambda)$ when $0\leq t\leq x$:
    \begin{equation}\label{annexe-resolvante-sl_est_tx}%
        |\G(x,t,\lambda)|\leq C
        \left(\frac{x}{1+|\omega|x}\right)^{a+1}\left(\frac{1+|\omega|t}{t}\right)^{a}\exp\left(|\im
        \omega|(x-t)\right),
    \end{equation}
    \begin{equation}\label{annexe-diffres-sl_est_tx}%
        \left|\frac{\partial\G}{\partial x}(x,t,\lambda)\right|\leq C
        \left(\frac{x}{1+|\omega|x}\right)^{a}\left(\frac{1+|\omega|t}{t}\right)^{a}\exp\left(|\im
        \omega|(x-t)\right).
    \end{equation}
    \item Estimations for the Green function $G(x,t,\lambda)$ when $0\leq x\leq t\leq
    1$:
    \begin{equation}\label{annexe-resolvante-sl_est_xt}%
        |\G(x,t,\lambda)|\leq C
        \left(\frac{1+|\omega|x}{x}\right)^{a}\left(\frac{t}{1+|\omega|t}\right)^{a+1}\exp\left(|\im
        \omega|(t-x)\right),
    \end{equation}
    \begin{equation}\label{annexe-diffres-sl_est_xt}%
        \left|\frac{\partial\G}{\partial x}(x,t,\lambda)\right|\leq
        C \left(\frac{1+|\omega|x}{x}\right)^{a+1}\left(\frac{t}{1+|\omega|t}\right)^{a+1} \exp\left(|\im
        \omega|(t-x)\right).
    \end{equation}
\end{itemize}

Rewriting relations $(54)-(56)$ in \cite{emot} section $7.2.8$
pp.$11-12$, we get
\begin{align}
    \label{annexe-diff-ja}
        x {j_a}'(x)=x j_{a-1}(x)-a j_a(x),\\
    \label{annexe-diff-ja-1}
    x {j_{a-1}}'(x)=a j_{a-1}(x)-x j_a(x),\\
    \label{annexe-diff-na}
    x \eta_{a}'(x)=x \eta_{a-1}(x)-a \eta_{a}(x),\\
    \label{annexe-diff-na-1}
    x \eta_{a-1}'(x)=a \eta_{a-1}(x)-x \eta_{a}(x).
\end{align}

We also deduce the following uniform estimates with $\omega\in\R,
|\omega|\rightarrow+\infty$
\begin{equation}\label{annexe-bessel-int-1}
    \displaystyle \int_0^1 j_a(\omega
 t)^2
 dt=\frac{1}{2}\left[1+\O{\frac{1}{\omega}}\right],
\end{equation}
\begin{equation}\label{annexe-bessel-int-2}
    \int_0^1 j_a(\omega t)\eta_a(\omega
t)dt=\O{\frac{1}{\omega}}.
\end{equation}

%   Bibliographie
%GATHER{Xbib.bib}   % For Gather Purpose Only
\bibliographystyle{plain}
\bibliography{Xbib}

\def\cprime{$'$}
\begin{thebibliography}{10}

\bibitem{bo}
G.~Borg.
\newblock Eine {U}mkehrung der {S}turm-{L}iouvilleschen {E}igenwertaufgabe.
  {B}estimmung der {D}ifferentialgleichung durch die {E}igenwerte.
\newblock {\em Acta Math.}, 78:1--96, 1946.

\bibitem{carlson2}
R.~Carlson.
\newblock Inverse spectral theory for some singular {S}turm-{L}iouville
  problems.
\newblock {\em J. Differential Equations}, 106(1):121--140, 1993.

\bibitem{carlson}
R.~Carlson.
\newblock A {B}org-{L}evinson theorem for {B}essel operators.
\newblock {\em Pacific J. Math.}, 177(1):1--26, 1997.

\bibitem{cs}
R.~Carlson and C.~Shubin.
\newblock Spectral rigidity for radial {S}chr\"odinger operators.
\newblock {\em J. Differential Equations}, 113(2):338--354, 1994.

\bibitem{emot}
A.~Erd{\'e}lyi, W.~Magnus, F.~Oberhettinger, and F.G. Tricomi.
\newblock {\em Higher transcendental functions. {V}ol. {II}}.
\newblock Robert E. Krieger Publishing Co. Inc., Melbourne, Fla., 1981.
\newblock Based on notes left by Harry Bateman, Reprint of the 1953 original.

\bibitem{gr}
J.-C. Guillot and J.~V. Ralston.
\newblock Inverse spectral theory for a singular {S}turm-{L}iouville operator
  on {$[0,1]$}.
\newblock {\em J. Differential Equations}, 76(2):353--373, 1988.

\bibitem{lev}
N.~Levinson.
\newblock The inverse {S}turm-{L}iouville problem.
\newblock {\em Mat. Tidsskr. B.}, 1949:25--30, 1949.

\bibitem{new}
R.~G. Newton.
\newblock {\em Scattering theory of waves and particles}.
\newblock Texts and Monographs in Physics. Springer-Verlag, New York, second
  edition, 1982.

\bibitem{ist}
J.~P{\"o}schel and E.~Trubowitz.
\newblock {\em Inverse spectral theory}, volume 130 of {\em Pure and Applied
  Mathematics}.
\newblock Academic Press Inc., Boston, MA, 1987.

\bibitem{rs}
M.~Reed and B.~Simon.
\newblock {\em Methods of modern mathematical physics. {II}. {F}ourier
  analysis, self-adjointness}.
\newblock Academic Press [Harcourt Brace Jovanovich Publishers], New York,
  1975.

\bibitem{rusa}
W.~Rundell and P.~E. Sacks.
\newblock Reconstruction of a radially symmetric potential from two spectral
  sequences.
\newblock {\em J. Math. Anal. Appl.}, 264(2):354--381, 2001.

\bibitem{zs}
L.~A. Zhornitskaya and V.~S. Serov.
\newblock Inverse eigenvalue problems for a singular {S}turm-{L}iouville
  operator on {$[0,1]$}.
\newblock {\em Inverse Problems}, 10(4):975--987, 1994.

\end{thebibliography}
\end{document}